\documentclass[final]{siamltex}
\usepackage{amsmath}
\usepackage{amssymb}
\usepackage{graphicx}
\usepackage{url}

\newcommand{\ttimes}{\!\times\!}

\title{Fifth Order Runge-Kutta-Nystr\"om Methods with Complex
	Coefficients\thanks{This work has been supported by a Veni
	Fellowship from Nederlandse Organisatie voor Wetenschappelijk
	Onderzoek (NWO) and by the Turkish Scientific Research Council
	(T\"UB\.ITAK) under Doktora Sonras\i\ Geri D\"on\"u\c{s}
	Program\i\ (2232)}}

\author{M. Atakan G\"urkan\thanks{Sterrewacht Leiden
    Niels Bohrweg 2, 2333 CA, Leiden, the Netherlands and
    Sabanci University, Faculty of Engineering and Natural Sciences
    Orhanl\i-Tuzla, 34956, \.Istanbul, Turkey
    ({\tt ato.gurkan@gmail.com}).}}

\begin{document}

\maketitle

\begin{abstract}
We present fifth order Runge-Kutta-Nystr\"om
methods, where we allow the timestep coefficients to assume complex values.
Among the methods with complex timesteps, we focus on the ones
with the coefficients that have positive real parts. This property
makes them suitable for problems where a negative coefficient is not
acceptable. In addition, the leading order terms in the error expansion
of these methods are purely imaginary, effectively increasing the
order of the methods by one for real variables.
\end{abstract}

\begin{keywords} 
Numerical integration, Runge-Kutta methods, Symplectic integrators,
Hamiltonian systems including symplectic integrators
\end{keywords}

\begin{AMS}
65D30,65L06,37M15,65P10
\end{AMS}

\pagestyle{myheadings}
\thispagestyle{plain}
\markboth{M. ATAKAN G\"URKAN}{COMPLEX RKN METHODS}

\section{Introduction}
\label{intro}
Splitting methods \cite{SEMA2008} provide a number of advantages for
studying the evolution of Hamiltonian systems. Not only are they simple
to implement, but also they can be fine tuned to exploit the structure
of the problem at hand \cite{2000AJ....120.2117L,2006gni.book.....H}. A large
class of physical problems are described by the separable Hamiltonian
\begin{equation}
\label{equ:hamiltonian}
H(q,p) = T(p) + V(q),
\end{equation}
where the kinetic energy $T(p)=\frac{1}{2} p^T M^{-1} p$ is quadratic
in momenta, and the potential energy $V(q)$ is only a function of
the coordinates.  For these problems, Hamilton's equations
\begin{equation}
\frac{d q_i}{d t} = \frac{\partial H}{\partial p_i},\quad
\frac{d p_i}{d t} = -\frac{\partial H}{\partial q_i}
\end{equation}
lead to the second order differential equation
\begin{equation}
\label{equ:basic}
\frac{d^2 q}{d t^2} = -M^{-1} \frac{\partial V(q)}{\partial q} = f(q)\,.
\end{equation}
This equation can be efficiently integrated
using Runge--Kutta--Nystr\"om (RKN) methods
\cite{1993sode.book.....H}.

RKN methods are particularly effective for orders higher
than 4, since some of the terms in the error expansion vanish
identically, thanks to the special structure in equation~\eqref{equ:basic}
\cite{1992MathComp....59.439O,Blanes_Moan}.
Zhu \& Qin \cite{WenJieZhu|MengZhaoQin199361} and Okunbor \& Skeel
\cite{Okunbor1994375} found
5th order explicit canonical RKN methods with 5 stages, which is the
minimum required by the order conditions.  The stages in these methods
have large (around unity) step size coefficients, some of which are negative . Large
coefficients can lead to large factors for integrations errors
\cite{Blanes_Moan}, while negative coefficients can make the
method unacceptable for the underlying problem \cite{Chin1997344}.
Starting from a splitting scheme, Chambers \cite{2003AJ....126.1119C} showed
that these difficulties can be overcome by allowing the step size
coefficients to be complex numbers, and found a third order splitting
method with two stages. The coefficients for this method have small
and positive real parts, leading to small errors. An interesting
property of this method is that the leading order terms in the error
expansion have purely imaginary coefficients. This makes the method
effectively one order higher for real variables (see ref. \cite{2010SEMA},
for other examples of such methods).

In this paper, we present multiple 5th order RKN methods with both
real and complex step size coefficients. We put special emphasis
on methods with coefficients that have small, positive real parts,
and purely imaginary leading error terms.

\section{5th Order Canonical RKN Methods}
\label{sec:methods}
For a system with Hamiltonian \eqref{equ:hamiltonian}, we can define
``velocity'' as $\dot{q} = M^{-1}p$. Then an $s$-stage RKN method is
given by \cite{1992MathComp....59.439O,Okunbor1994375}
\begin{equation}
\begin{split}
y_i = q_n + c_ih\dot{q}_n + h^2 \sum_{j=1}^s a_{ij} f(y_j), \quad
	i=1,2,\ldots,s\,,\\
q_{n+1} = q_n + h\dot{q}_n + h^2\sum_{i=1}^s b_i f(y_i),\quad
\dot{q}_{n+1} = \dot{q}_n + h\sum_{i=1}^s B_i f(y_i)\,,
\end{split}
\end{equation}
where $f(q)$ is defined in equation \eqref{equ:basic}. This method is
explicit, that is a step depends only on previous steps, if $a_{ij}=0$
for $j\ge i$. An explicit $s$-stage RKN method without redundant stages
is canonical if \cite{1992MathComp....59.439O,Okunbor1994375}
\begin{align}
b_i = B_i(1-c_i),        \quad 1&\le i \le s,\\
a_{ij} = B_j(c_i - c_j), \quad j&<i\,.
\end{align}
For methods of order $5$, the following order conditions must hold
\cite{Okunbor1994375,WenJieZhu|MengZhaoQin199361}:
\begin{align}
\label{equ:conditions}
&t_1:\; \sum_i B_i = 1, 
    & t_6:&\; \sum_i \sum_{j<i} B_i B_j (c_i -c_j) =
	\frac{1}{6},\notag\\
&t_2:\; \sum_i B_i c_i = \frac{1}{2},
    &t_7:&\; \sum_i \sum_{j<i} B_i B_j c_i (c_i-c_j) =
	\frac{1}{8},\notag\\
&t_3:\; \sum_i B_i c_i^2 = \frac{1}{3},
    &t_8:&\; \sum_i \sum_{j<i} B_i B_j c_i^2(c_i-c_j) =
	\frac{1}{10},\\
&t_4:\; \sum_i B_i c_i^3 = \frac{1}{4},
    &t_{9}:&\; \sum_i\sum_{j<i} B_i B_j c_i c_j (c_i-c_j) =
	\frac{1}{30},\notag\\
&t_5:\; \sum_i B_i c_i^4 = \frac{1}{5},
    &t_{10}:&\; \sum_i \sum_{j<i} \sum_{l<i} B_i B_j B_l (c_i-c_j)
                            (c_i-c_l) = \frac{1}{20}\,.\notag
\end{align}

By solving these equations, we can obtain an RKN method of order $5$;
note, however, that these equations by themselves do not give us the
error behaviour.

\section{Splitting Scheme}
Since we will be using a Hamiltonian splitting scheme for
implementing our integrators and error analysis, we now briefly review
the basics of this technique using Lie algebra, in a manner
similar to the treatment of Yoshida \cite{1993CeMDA..56...27Y}.

First we write Hamilton's equations for a state $w(t)\equiv(p(t),q(t))$ as
\begin{equation}
\frac{d w}{d t}= \{w,H\} = \left(\sum_i
    \frac{\partial H}{\partial p_i} \frac{\partial\:}{\partial q_i}
   -\frac{\partial H}{\partial q_i} \frac{\partial\:}{\partial p_i}
\right) w = \mathbb{H} w = (\mathbb{T} + \mathbb{V})w\,,
\end{equation}
where $\mathbb{H}$, $\mathbb{T}$, $\mathbb{V}$ are Lie operators
\cite{Chin1997344} corresponding to the Hamiltonian, the kinetic
energy and the potential energy, respectively. The solution for $w(t)$
can then be written as
\begin{equation}
\label{equ:formal}
w(\tau) = e^{\tau\mathbb{H}} w(0)\,,
\end{equation}
and the exponential operator is approximated by
\begin{equation}
\label{equ:expand}
e^{\tau\mathbb{H}} \approx e^{\tau\mathbb{Z}} = \prod_i e^{\alpha_i \tau \mathbb{T}} 
                             e^{\beta_i \tau \mathbb{V}}\,.
\end{equation}
The expansion for $\mathbb{Z}$ in terms of $\mathbb{T}$
and $\mathbb{V}$ can be obtained by repeatedly applying the
Baker-Campbell-Hausdorff (BCH) formula (for BCH expressions with
minimum number of terms, see \cite{Blanes2004135,2009JMP....50c3513C}).
The coefficients $\alpha_i$, $\beta_i$ are chosen such that
$\tau\mathbb{Z} = \tau\mathbb{H}+ O(\tau^{n+1})$, for a method of
order $n$. Note that the expansion of $\mathbb{Z}$ can be considerably
simplified, since the commutator
$[\mathbb{V}, [\mathbb{V}, [\mathbb{T}, \mathbb{V}]]]$ vanishes for the
systems described by equation~\eqref{equ:basic}.

The evolution due to the exponential
operators on the right hand side of equation~\eqref{equ:expand} is given by
simple displacements
\begin{equation}
e^{\alpha_i \tau \mathbb{T}} (q,p) = (q+\alpha_i \tau p, p)
\quad\mathrm{and}\quad
e^{\beta_i \tau \mathbb{V}} (q,p) = (q, p +\beta_i \tau f(q))\,.
\end{equation}
The relation between the coefficients of
an RKN method and a splitting scheme is given by  \cite{1992MathComp....59.439O}:
\begin{equation}
\alpha_i = c_i - c_{i-1},\: \alpha_{s+1} = 1-c_s, \beta_i = B_i,\:
\beta_{s+1} = 0\,,
\end{equation}
(where $1\le i \le s$, $c_0=0$),
leading to the following procedure for the method \cite{1992MathComp....59.439O}:
\begin{equation*}
\begin{split}
&y_0 = q_n,\\
&\dot{y}_0 = \dot{q}_n,\\
&\mathsf{for}\; i = 1, 2,\ldots,s\\
&\quad y_i = y_{i-1} + h(c_i - c_{i-1})\dot{y}_{i-1},\qquad\qquad\qquad\qquad\qquad\qquad\qquad\qquad\\
&\quad \dot{y}_i = \dot{y}_{i-1} + h B_i f(y_i),\\
&q_{n+1} = y_s + h(1-c_s)\dot{y}_s\\
&\dot{q}_{n+1} = \dot{y}_s\,.
\end{split}
\end{equation*}

This corresponds to a splitting scheme with the structure
\begin{equation}
\label{equ:rkna}
e^{\tau\mathbb{H}} \approx e^{\tau\mathbb{Z}} =
e^{\alpha_6 \tau \mathbb{T}} e^{\beta_5 \tau \mathbb{V}}
 e^{\alpha_5 \tau \mathbb{T}} \ldots e^{\beta_1 \tau \mathbb{V}}
 e^{\alpha_1 \tau \mathbb{T}}\,,
\end{equation}
to which we will refer as \texttt{RKNA} scheme.
Another method with the structure
\begin{equation}
\label{equ:rknb}
e^{\tau\mathbb{H}} \approx e^{\tau\mathbb{Z}} =
 e^{\beta_6 \tau \mathbb{V}} e^{\alpha_5 \tau \mathbb{T}} 
 e^{\beta_5 \tau \mathbb{V}} \ldots e^{\alpha_1 \tau \mathbb{T}} 
 e^{\beta_1 \tau \mathbb{V}} \,,
\end{equation}
to which we will refer as \texttt{RKNB} scheme, would have similar
computational cost. However, because of the special
structure of the RKN methods, the coefficients for the latter scheme would be
entirely different \cite{Blanes_Moan}. To obtain the coefficients from
the same order conditions, we
rewrite equation \eqref{equ:rknb} as
\begin{equation}
e^{\tau\mathbb{H}} \approx e^{\tau\mathbb{Z}} =
 e^{\alpha_7 \tau \mathbb{T}} 
 e^{\beta_6 \tau \mathbb{V}} e^{\alpha_6 \tau \mathbb{T}} 
 e^{\beta_5 \tau \mathbb{V}} \ldots e^{\alpha_2 \tau \mathbb{T}} 
 e^{\beta_1 \tau \mathbb{V}} e^{\alpha_1 \tau \mathbb{T}} \,,
\end{equation}
i.e., as a 6 stage method. Since $\alpha_7 = \alpha_1 = 0$, we have
$c_1 = 0$ and $c_6=1$. Hence, the number of unknowns is again equal to
the number of equations.

\section{Constructing the Methods}
We solved the order conditions, eq. \eqref{equ:conditions},
using {\tt fsolve} routine of MAPLE, with {\tt complex}
optional keyword. We started from a large number of random initial
guesses for $B_i$ and $c_i$ in the rectangle $-2 \le \mathfrak{Re}(z),
\mathfrak{Im}(z) \le 2$. This yielded both of the previously known
real solutions for \texttt{RKNA} scheme, along with their adjoints;
as well as three real solutions for \texttt{RKNB} scheme that seems to
have not been discovered before. In addition, we found a large number of
complex solutions for both schemes.

Once the coefficients are calculated, we repeatedly applied BCH
formula to calculate the leading order error terms as nested commutators of
$\mathbb{T}$ and $\mathbb{V}$.  Because of the Jacobi identity and the
simplification $[\mathbb{V}, [\mathbb{V}, [\mathbb{T}, \mathbb{V}]]]
= 0$, not all commutators are independent. Also, to calculate the
leading order error of the method, we are only interested
in terms with up to 6 operators; so, we worked in a Philip Hall basis
\cite{2009JMP....50c3513C}
with 2 generators and nilpotency 7, leading to 23 terms. Denoting
our generators by $X_1 = \mathbb{T}$ and $X_2=\mathbb{V}$, we 
constructed a ``multiplication table'':
\begin{equation}
\label{equ:mult_table}
\begin{split}
&[X_1, X_2]    = X_3;\,
[X_1, X_3]    = X_4;\,
[X_1, X_4]    = X_6;\,
[X_1, X_5]    = X_7;\,\\ &
[X_1, X_6]    = X_{12};\,
[X_1, X_7]    = X_9 + X_{13};\,
[X_1, X_9]    = X_{16};\, \\ &
[X_1, X_{10}] = X_{11} + X_{17};\,
[X_1, X_{12}] = X_{19};\,
[X_1, X_{13}] = X_{16} + X_{20};\, \\ & 
[X_1, X_{14}] = -X_{11} - X_{17};\,
[X_2, X_3]    = X_5;\,
[X_2, X_4]    = X_7;\,
[X_2, X_6]    = X_{13};\, \\ &
[X_2, X_7]    = -X_{10};\, 
[X_2, X_9]    = -X_{11} + X_{17};\,
[X_2, X_{12}] = X_{20};\, \\ &
[X_2, X_{13}] = -3X_{17};\,
[X_3, X_4]    = X_9;\,
[X_3, X_5]    = X_{10};\, 
[X_3, X_6]    = X_{16};\, \\ &
[X_3, X_7]    = X_{17};\,
[X_4, X_5]    = X_{11};
\quad \mathrm{with}\quad [X_i, X_j] = -[X_j, X_i]\,.
\end{split}
\end{equation}
Apart from the commutators given in this table, all commutators of
the basis operators vanish. Using this table, we wrote a
{\tt Python}\footnote{\protect\url{http://www.python.org}} program, using the
{\tt SymPy}\footnote{\protect\url{http://code.google.com/p/sympy}} package, to
obtain the expansion for $\mathbb{Z}$ up to and including
6 operator terms. This allowed us to validate the methods and
calculate the coefficients for the leading order error.

For \texttt{RKNA} scheme, we found two previously published methods
with real coefficients. For \texttt{RKNB} scheme, we found three
methods with real coefficients that we did not find elsewhere in
the literature. We give the coefficients for these methods in Table
\ref{tab:real_coeffs}.

We also found a large number of methods with complex coefficients. Among
these, ten for \texttt{RKNA} scheme and five for \texttt{RKNB}
scheme had coefficients with positive real parts and purely imaginary
leading order errors. We give the coefficients for two methods with
smallest error coefficients for each scheme in Table
\ref{tab:comp_coeffs}. 

The leading order errors are given by the coefficients in front of
the six-term commutators in the expansion of $\mathbb{Z}$. These
coefficients, for the methods presented, are given in Table
\ref{tab:error_terms}.

\begin{table}
\renewcommand{\arraystretch}{0.9}
\caption{Coefficients for real fifth order splitting schemes. The
    adjoint methods are also valid and can be obtained 
    by reversing the order of the coefficients.
    The coefficients are not entirely
    independent since $\sum \alpha_i = \sum \beta_i = 1$.}
\label{tab:real_coeffs}
\begin{tabular}{lrrrr}
 & \multicolumn{1}{c}{$\alpha_i$} & \multicolumn{1}{c}{$\beta_i$} \\
\hline
\texttt{AR1} &  0.96172990014645096  &   0.39682804502722538  \\
             & -0.09525408032034999  &  -0.824377563589592    \\
             & -0.73942683539212613  &   0.2042028689314904   \\
             &  0.62730935078241887  &   1.0021847152077973   \\
             & -0.52506178465602220  &   0.22116193442307898  \\
             &  0.77070344943962849  &   -- \\
\texttt{AR2} &  0.69883375727545265  &   0.40090379269659899  \\
             & -0.49469565362085154  &   0.95997088013405985  \\
             &  0.81641946634957295  &   0.0884951581272243   \\
             & -0.65762956677338285  &   1.2214390923487315   \\
             & -0.057841894299102682 &  -1.6708089233066146   \\
             &  0.69491389106831146  &   -- \\
\texttt{BR1} &  0.54200976680171613  &   0.24566294009066009  \\
             & -0.04060817665564392  &   1.1433587581365421   \\
             & -0.87779698530109766  &  -1.3796706973507000   \\
             &  0.86474236062251646  &  -0.019611260781217307 \\
             &  0.51165303453250898  &   0.87087215441178844  \\
             &  --                   &   0.13938810549292669  \\
\texttt{BR2} &  0.42637413177222316  &   0.15102308452230116  \\
             & -0.82438794434938248  &   0.72768821316253478  \\
             & -0.63140077574154094  &  -0.26217627934521390  \\
             &  0.38590710518893978  &  -0.044211509719803855 \\
             &  1.6435074831297605   &   0.23596222045571453  \\
             &  --                   &   0.19171427092446728  \\
\texttt{BR3} &  1.0413749845202060   &   0.12696076271851077  \\
                & -0.61784769849171965  &  -1.4166626058695677   \\
                &  0.62570540985789957  &  -0.62172666654176438  \\
                & -0.63446409452971410  &   0.69301448863793809  \\
                &  0.58523139864332822  &   1.2079876026916669   \\
                &  --                   &   1.0104264183632164   \\
\end{tabular}
\end{table}
\begin{table}
\caption{Coefficients for a selection of complex fifth order splitting
    schemes, with purely imaginary leading order error terms. All the
    methods are skew-symmetric, i.e., reversing the order of the
    coefficients gives their complex conjugates. This property and
    $\sum_i \alpha_i = \sum_i \beta=1$ can be used to calculate the
    missing coefficients. The adjoints of the
    methods, obtained by reversing the order of the coefficients,
    are also valid.
}
\label{tab:comp_coeffs}
\begin{tabular}{rlrr}
 & & \multicolumn{1}{c}{$\mathfrak{Re}$} &
     \multicolumn{1}{c}{$\mathfrak{Im}$}   \\
\hline
\texttt{AC1}& $\alpha_1$ & 0.087808410045663212 &  0.028523844251341822 \\
            & $\alpha_2$ & 0.17916539354193987  & -0.067857083007249973 \\
            & $\alpha_3$ & 0.23302619641239692  & -0.097952003128893425 \\
            & $\beta_1$  & 0.17526734338348050  &  0.057642040076250593 \\
            & $\beta_2$  & 0.18488007701471166  & -0.19410647329733509  \\
\texttt{AC2}& $\alpha_1$ & 0.087634204536037057 &  0.028807372065269351 \\
            & $\alpha_2$ & 0.18007104463252914  & -0.068253589313355443 \\
            & $\alpha_3$ & 0.23229475083143381  & -0.097060961378624794 \\
            & $\beta_1$  & 0.17526840907207411  &  0.057614744130538702 \\
            & $\beta_2$  & 0.18487368019298416  & -0.19412192275724959  \\
\texttt{BC1}& $\beta_1$  & 0.093106790861751605 & -0.026812950639104607 \\
            & $\beta_2$  & 0.14578332225686154  &  0.076033669531385746 \\
            & $\beta_3$  & 0.26110988688138685  &  0.10851236434561279  \\
            & $\alpha_1$ & 0.15950063058390336  & -0.060127448366782494 \\
            & $\alpha_2$ & 0.19085044206705213  &  0.20369642527600502  \\
\texttt{BC2}& $\beta_1$  & 0.10625796854753310  & -0.037213537431233983 \\
            & $\beta_2$  & 0.35767992721948460  & -0.022169204268009056 \\
            & $\beta_3$  & 0.036062104232982296 &  0.057072185585748646 \\
            & $\alpha_1$ & 0.26934942679787788  & -0.093675141997563700 \\
            & $\alpha_2$ & 0.14580813747862993  &  0.49930185549019606  \\
\end{tabular}
\end{table}
\begin{table}
\caption{Coefficients of the terms in the expansion 
    $\mathbb{H} - \mathbb{Z}$, for the methods given in
    Tables \ref{tab:real_coeffs} and \ref{tab:comp_coeffs}.
}
\label{tab:error_terms}
\begin{tabular}{lrrrrrr}
 & \multicolumn{1}{c}{$X_{11}$}
 & \multicolumn{1}{c}{$X_{16}$} 
 & \multicolumn{1}{c}{$X_{17}$} 
 & \multicolumn{1}{c}{$X_{19}$} 
 & \multicolumn{1}{c}{$X_{20}$} 
 & \multicolumn{1}{c}{$\left(\sum_i |X_i|^2\right)^{1/2}$} \\
\hline
{\tt AR1} & $-3.1\ttimes10^{-3}$ & $ 4.5\ttimes10^{-4}$
          & $ 5.6\ttimes10^{-3}$ & $ 1.8\ttimes10^{-5}$ 
          & $-6.3\ttimes10^{-5}$ & $6.44\ttimes 10^{-3}$ \\
{\tt AR2} & $ 5.0\ttimes10^{-3}$ & $-3.1\ttimes10^{-4}$
          & $-5.7\ttimes10^{-3}$ & $ 1.5\ttimes10^{-5}$ 
          & $ 2.5\ttimes10^{-4}$ & $7.58\ttimes 10^{-3}$ \\
{\tt BR1} & $ 1.5\ttimes10^{-4}$ & $ 6.2\ttimes10^{-5}$
          & $-4.1\ttimes10^{-4}$ & $ 5.7\ttimes10^{-5}$
          & $ 1.6\ttimes10^{-4}$ & $ 4.71\ttimes10^{-4}$ \\
{\tt BR2} & $-7.0\ttimes10^{-5}$ & $ 9.7\ttimes10^{-4}$
          & $ 2.7\ttimes10^{-4}$ & $ 5.2\ttimes10^{-4}$
          & $ 5.2\ttimes10^{-4}$ & $1.25\ttimes10^{-3}$ \\
{\tt BR3} & $ 3.5\ttimes10^{-2}$ & $ 1.9\ttimes10^{-4}$
          & $-3.5\ttimes10^{-2}$ & $-3.0\ttimes10^{-5}$ 
          & $-1.0\ttimes10^{-4}$ & $4.90\ttimes10^{-2}$ \\
{\tt AC1} & $-i2.5\ttimes10^{-6}$ & $ i7.3\ttimes10^{-6}$
          & $-i5.0\ttimes10^{-6}$ & $ i8.3\ttimes10^{-7}$
          & $ i4.2\ttimes10^{-6}$ & $1.02\ttimes10^{-5}$ \\
{\tt AC2} & $ i3.0\ttimes10^{-6}$ & $-i7.9\ttimes10^{-6}$
          & $ i6.0\ttimes10^{-6}$ & $-i9.4\ttimes10^{-7}$
          & $-i4.7\ttimes10^{-6}$ & $1.15\ttimes10^{-5}$ \\
{\tt BC1} & $ i1.2\ttimes10^{-5}$ & $-i1.3\ttimes10^{-5}$
          & $ i8.4\ttimes10^{-6}$ & $ i2.2\ttimes10^{-7}$ 
          & $-i4.2\ttimes10^{-6}$ & $1.96\ttimes10^{-5}$ \\
{\tt BC2} & $-i4.1\ttimes10^{-5}$ & $ i4.9\ttimes10^{-5}$
             & $-i1.2\ttimes10^{-4}$ & $ i7.7\ttimes10^{-6}$ 
	     & $ i5.0\ttimes10^{-5}$ & $1.49\ttimes10^{-4}$ \\
\end{tabular}
\end{table}

\section{Numerical Experiments}
As a simple validation, we first compare the behaviour of two of these
methods (\texttt{RKNAC1} and \texttt{RKNBR1}) with
other methods from the literature, for the gravitational two-body
problem. We integrated the equations of motion of two equal
point masses, on orbit around each other with eccentricity $e=0.2$,
for fifty orbital periods. To follow the orbit, we used a
Gragg-Bulirsch-Stoer (GBS) integrator \cite{1993sode.book.....H}.
At each timestep, we also calculated the
expected position and velocity for each of the methods we tested. The
difference between the outcomes of GBS integrator and the method gives
an estimate of the error made, as long as they are not dominated by
the truncation errors arising from limited machine accuracy
($\sim 10^{-15}$). We then calculated $0.25$ and $0.75$ quantile and took
their difference, to get the {\em interquartile range}. This is a robust
statistic that gives a measure of the dispersion in data.

\begin{figure}
\includegraphics[angle=270,width=\textwidth,clip]{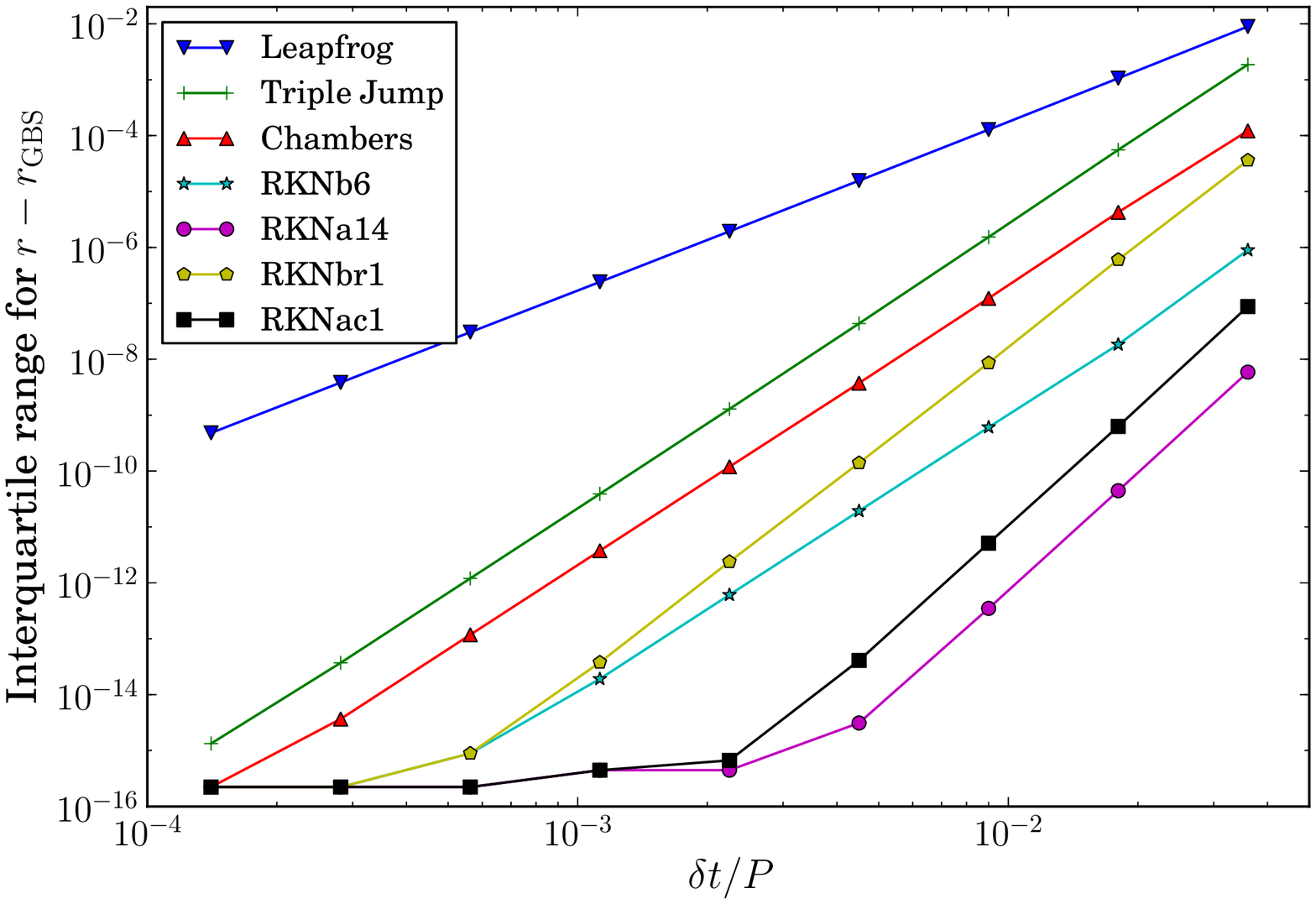}
\caption{Comparison of the behaviour of various methods with respect
to changing step size, for gravitational two-body problem. $x$-axis
the logarithm of the time step measured in orbital periods, $y$-axis is
the logarithm of the interquartile range for position error.  {\tt Leapfrog} 
(down-triangles) is the standard second order St\"ormer-Verlet
scheme, {\tt Triple jump}
(pluses) is the fourth order scheme of Yoshida \cite{1993CeMDA..56...27Y},
{\tt Chambers} (up-triangles) is the complex third order scheme of Chambers
\cite{2003AJ....126.1119C},{\tt RKNb6} (stars) and {\tt RKNa14} 
(circles) are optimized fourth and sixth order methods by Blanes \& Moan 
\cite{BlanesMoan2001,Blanes_Moan}, {\tt RKNbr1} (pentagons) and {\tt
RKNac1} (squares) are two of the methods presented in this work.
}
\label{fig:order}
\end{figure}

In figure~\ref{fig:order}, we plot the interquartile range
of position error, for various step sizes and different methods.
For implementing the method of Chambers \cite{2003AJ....126.1119C} and
our methods,
we chose to throw away the imaginary part of the positions and the
velocities after each step. This destroys the symplecticness of
the methods (or rather reduces the order for which the methods are
symplectic) but leads to good error behaviour \cite{2010SEMA}.

The comparison indicates that, for this problem, method of Chambers
\cite{2003AJ....126.1119C} shows fourth order and method \texttt{AC1}
shows sixth order behaviour, even though they are formally third
and fifth order, respectively. Analogous results were obtained by
Chambers \cite{2003AJ....126.1119C} for similar problems. It is interesting
to see that the optimized 4th order method \texttt{RKNb6}
\cite{Blanes_Moan} outperforms a method of higher order, \texttt{RKNbr1},
down to machine accuracy level.

\begin{figure}
\includegraphics[angle=270,width=\textwidth,clip]{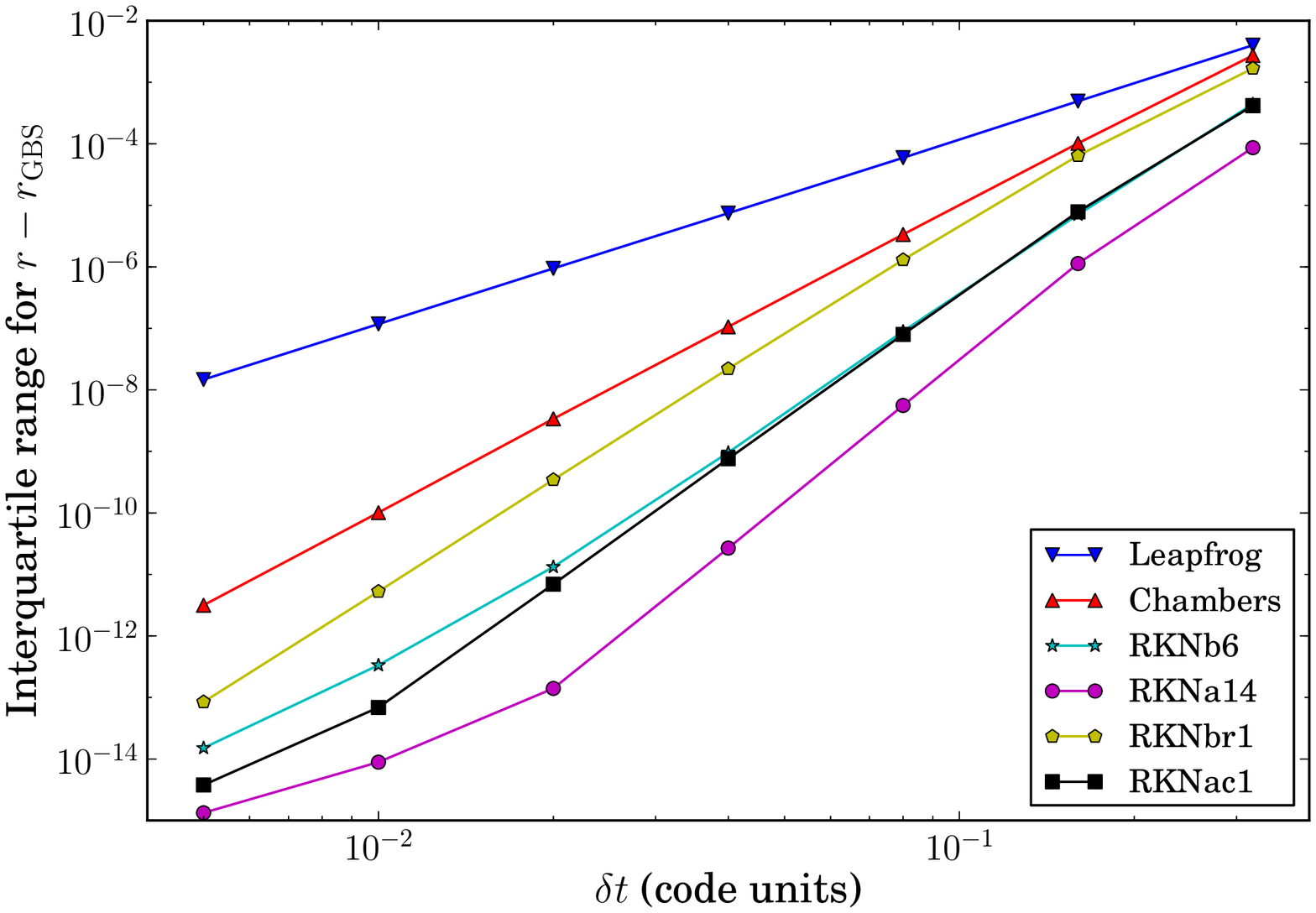}
\caption{Comparison of the behaviour of various methods with respect
to changing step size, for gravitational 400-body problem. $x$-axis
the logarithm of the time step in code units, $y$-axis is
the logarithm of the interquartile range for position error. The symbols
correspond to the same methods as figure \ref{fig:order}, except
\texttt{Triple Jump} is omitted, since it nearly coincides with
\texttt{Chambers} for this problem.}
\label{fig:order400}
\end{figure}

To make a more comprehensive and meaningful test, we also developed 
integrators for the  gravitational $N$-body problem, based on various methods.
We simulated the evolution of a Plummer sphere with 400 equal mass
particles (see ref. \cite{1974A&A....37..183A} for a procedure for constructing a
Plummer sphere). We followed the same procedure to estimate the errors
and calculated the interquartile range for various methods and
stepsizes.  The dependence of errors on step size is given in figure
\ref{fig:order400}. 

Since complex arithmetic requires more operations than real arithmetic,
we also made a comparison of CPU times for various methods. It was not
possible to calculate CPU times for two-body or 400-body problems
accurately, since
each step took too little time. Consequently, we setup a system of 10000
particles and integrated over a few steps. While calculating the CPU
time, we subtracted off the time spent for the last substep for an integration
step, since all the schemes we consider have the so called first-same-as-last
property. For example, in \texttt{RKNA} schemes (equation
\ref{equ:rkna}) $e^{\alpha_1 \tau \mathbb{T}}$ substep can be combined
with $e^{\alpha_6 \tau \mathbb{T}}$ substep of the next step. We present
the CPU times per step for various methods in table \ref{tab:CPU_times}.

\begin{table}
\caption{CPU times spent at each step for various methods. Values are
    normalized to Leapfrog method's CPU time.}
\label{tab:CPU_times}
\begin{tabular}{lr}
\texttt{Leapfrog} & 1 \\
\texttt{Chambers} & 12  \\
\texttt{Triple Jump} & 3  \\
\texttt{RKNb6} & 6  \\
\texttt{RKNa14} & 14  \\
\texttt{RKNbr1} & 5  \\
\texttt{RKNac1} & 28  \\
\end{tabular}
\end{table}

\begin{figure}
\includegraphics[angle=270,width=\textwidth,clip]{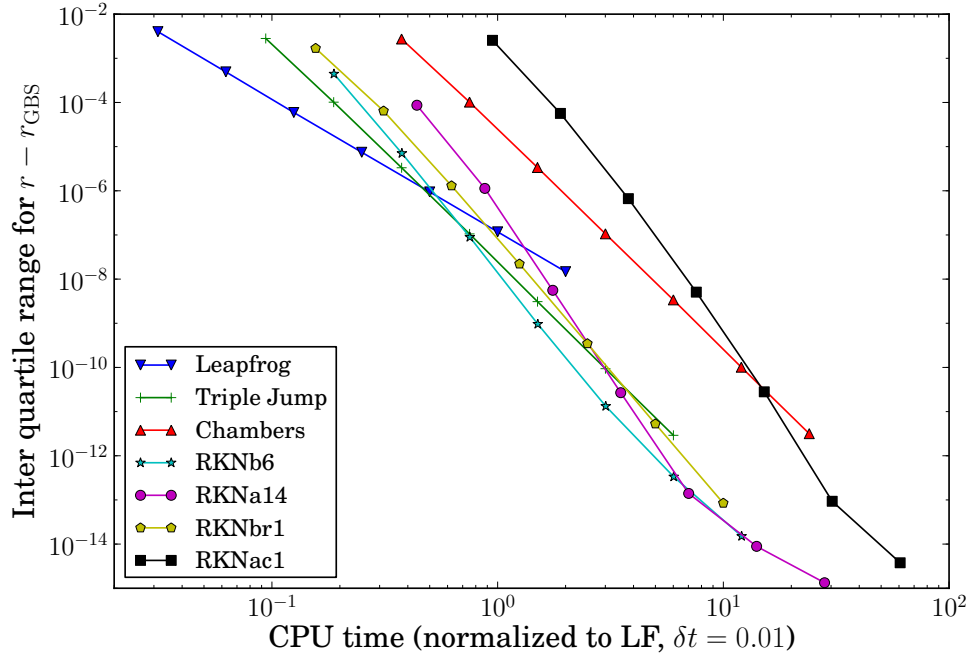}
\caption{CPU time spent vs. error for various methods applied to
gravitational 400-body problem. $x$-axis
the logarithm of the CPU time, normalized to time spent for Leapfrog
method at time step $\delta t = 0.01$, $y$-axis is
the logarithm of the interquartile range for position error. The symbols
correspond to the same methods as figure \ref{fig:order}.}
\label{fig:cputime}
\end{figure}

The data here show that the integration time is proportional to the
number of stages and using complex arithmetic increases the computation
cost by about a factor of 6. This last factor surely depends on the
problem, the implementation and the compiler. We present the error vs.
CPU time, based on these findings in figure \ref{fig:cputime}.


\section{Discussion}
In this work we constructed fifth order RKN methods with complex
coefficients. Some of the methods we found have much smaller timesteps
than the previously known methods, leading to smaller leading order
errors. In addition, many of them have positive real
coefficients, making them suitable for problems
where negative real timesteps are not acceptable
\cite{PhysRevA.42.6991,Bandrauk1991428,Chin1997344}. Similar methods,
satisfying this requirement, were also developed by Bandrauk et
al. \cite{Bandrauk2006346}. Other high order methods with complex
coefficients with small positive real parts were developed and analyzed by 
Hansen \& Ostermann \cite{MR2545819} and Castella et al. \cite{MR2545817};
however these authors did not specifically study the RKN case, which
allows considerable simplification for high order schemes.

For problems where complex arithmetic and/or positive real parts of the
coefficients is a necessity, we expect these methods to be
already competitive with lower order methods. However, the results of
Blanes \& Moan \cite{Blanes_Moan} and Sofroniou \& Spaletta
\cite{2005Concat..S}
suggest that there is room for improvement by increasing the number of
stages. We consider finding fully optimized methods with more stages
beyond the scope of this paper; such an effort would require developing
new software and would need spending considerable computational power.
However, the methods found here can be improved by readily
available tools. The idea is to turn a 5-stage method into a 6-stage one
by adding an additional stage. 
\begin{equation}
e^{\alpha_6 \tau \mathbb{T}} e^{\beta_5 \tau \mathbb{V}}
 e^{\alpha_5 \tau \mathbb{T}} \ldots e^{\beta_1 \tau \mathbb{V}}
 e^{\alpha_1 \tau \mathbb{T}}
 \rightarrow\\
 e^{\beta'_7 \tau \mathbb{V}} e^{\alpha_6 \tau \mathbb{T}} 
 e^{\beta'_6 \tau \mathbb{V}} \ldots e^{\alpha_1 \tau \mathbb{T}} 
 e^{\beta'_1 \tau \mathbb{V}}\,,
\end{equation}
and set $\beta'_7 = \beta'_1 = 0$. We can then start in the vicinity
of this solution and use MAPLE's minimization routines. Starting from
solution \texttt{AC1} in table \ref{tab:comp_coeffs}, we obtained the following
skew-symmetric method:
\begin{equation}
\begin{split}
 \alpha_1 & = .101907705405177865 + i.130701756906677735\,,\\ 
 \alpha_2 & = .218628781976265590 + i.0126440811480678494\,,\\
 \alpha_3 & = .179463512618556560 - i.148112326926992222\,, \\
 \beta'_1 & = .0489489561074426954 + i.0669384556781967844\,,\\
 \beta'_2 & = .166479171860817010  + i.0764027877516731402\,,\\
 \beta'_3 & = .192297943665939275  - i.0835834606213808479\,,\\
 \beta'_4 & = .184547856731601789 \,.
\end{split}
\end{equation}
This method has an error $\left(\sum_i |X_i|^2\right)^{1/2} = 2.6\times
10^{-7}$, about two orders of magnitude smaller than the other methods
(see table \ref{tab:error_terms}).

Since we minimized the (imaginary) sixth order error terms, this
optimization does not benefit the solutions of gravitational $N$-body
problem. A possible venue for exploration is to increase the number of
stages and minimize the (real part of) seventh order errors by using
the extra variables.  For this investigation, using a Lyndon basis
would be more preferable, since the BCH expansion has much fewer terms
in this basis \cite{2009JMP....50c3513C} and the consequences of the
simplification $[\mathbb{V}, [\mathbb{V}, [\mathbb{T}, \mathbb{V}]]]
= 0$ is more straightforward. Here, we used a Philip Hall basis,
since this allowed us to check our algebra by the ``Lie Tools
Package''\footnote{\protect\url{http://www.cim.mcgill.ca/~migueltt/ltp/ltp.html}}
of Miguel Torres-Torriti.

The source code of the \texttt{MAPLE}, \texttt{Python} and \texttt{C}
programs used in this work are freely available online:\
\url{https://github/atakan/Complex_Coeff_RKN}

\section*{Acknowledgments}
I thank Yuri Levin, Sergio Blanes, Ander Murua, Fernando Casas, Tevhide
Altekin and Tolga G\"uver for various discussions and an anonymous referee
for suggestions that improved this paper.

\bibliographystyle{siam}
\bibliography{RKN5}  

\end{document}